\documentclass[11pt]{amsart}

\usepackage{amsfonts}
\usepackage{amscd}
\usepackage{amssymb}
\usepackage{amsmath}
\usepackage{latexsym} 
\usepackage{graphicx}
\usepackage{color}
\usepackage{epsfig}
\usepackage{url}

\newtheorem{prop}{Proposition}[section]

\newtheorem{rem}[prop]{Remark}
\newtheorem{ex}[prop]{Example}
\newtheorem{cor}[prop]{Corollary}
\newtheorem{thm}[prop]{Theorem}

\def\KK{{\mathbb{K}}} \def\PP{{\mathbb{P}}} \def\RR{{\mathbb{R}}}
 \def\NN{{\mathbb{N}}} \def\ZZ{{\mathbb{Z}}}

\def\Kc{{\mathcal{K}}} \def\Cc{{\mathcal{C}}} \def\Dc{{\mathcal{D}}}
  \def\OO{{\mathcal{O}}}
  \def\MM{{\mathcal{M}}}
  
  \def\EE{{\mathcal{E}}}

\def\fg{\mathbf{f}}  
\def\xg{\mathbf{x}} \def\dg{\mathbf{d}} \def\kg{\mathbf{k}}

\def\mg{\mathbf{m}}

\def\mm{\mathfrak{m}}

\def\Res{\mathrm{Res}}

\def\dim{\mathrm{dim}}
\def\det{\mathrm{det}}
\def\depth{\mathrm{depth}}
\def\deg{\mathrm{deg}}

\def\Spec{\mathrm{Spec}}

\def\rank{\mathrm{rank}}

\def\Im{\mathrm{Im}}

\title{A resultant approach to detect
intersecting curves in $\PP^3$}
\date{\today}
\author{Laurent Bus\'e and Andr\'e Galligo}
\address{
  Universit\'e de Nice - Sophia Antipolis,\newline Laboratoire de
  Math\'ematiques, Parc Valrose 06108 Nice Cedex 2.}
\email{{lbuse,galligo}@unice.fr}
\thanks{This work is partially supported by the european projects GAIA
  II IST-2001-35512, and ECG IST-2000-26473.}

\begin{document}

\begin{abstract}
  Given two curves in $\PP^3$, either implicitly or by a
  parameterization, we want to check if they intersect. For that
  purpose, we present and further develop generalized resultant
  techniques. Our aim is to provide a closed formula in the inputs
  which vanishes if and only if the two curves intersect. This could
  be useful in Computer Aided Design, for computing the intersection
  of algebraic surfaces.
\end{abstract}

\maketitle

\section{Introduction}
Resultants have a long history going back to B\'ezout in the
eighteenth century. This algebraic construction received a renewed
interest in the last decade with many new contributions and
extensions. Generalized resultant is now a major tool in elimination
theory, both for theoretical and practical purposes
\cite{GKZ94,CLO98}.

In the homogeneous case, a generalized resultant is a closed formula
which answer a problem of the following kind~: let $X$ be an
irreducible projective variety in $\PP^n$,
$\underline{a}=(a_1,\ldots,a_N)$ be parameters, and
$f_1(x,\underline{a})$,\ldots,$f_s(x,\underline{a})$ be $s$
homogeneous polynomials in $x\in X$. A resultant, if it exists, is a
polynomial in $\underline{a}$, called $\Res(a)$, such that
$$\Res(a)=0 \Leftrightarrow \{ x\in X \ / \ f_i(x,\underline{a})=0
\}\neq \emptyset.$$
Most of the works in the subject focus on the case
$s=\dim(X)+1$. However the first author studied in its PhD thesis
\cite{BusPhD} (see also \cite{Bus02}) \emph{determinantal} resultants
which can be used for some cases with $s > \dim(X)+1$. Motivated by
applications in CAD (Computer Aided Design) we now go further in this
direction.

Indeed, a major problem in CAD is the efficient computation of the
intersection of two algebraic surface patches in $\RR^3$ given by
parameterizations. A popular strategy is to first determine the
implicit representation of one of the two surfaces, therefore
implicitization techniques based on generalized resultants have been
developed, including sparse resultants. Unfortunately, they bring up
to the difficult question of the description of the base points of the
parameterization \cite{Bus01,Cox01,BuJo02}.

The second author proposed in \cite{Gal03} to compute only a
semi-implicit representation of a parameterized surface, which is much
easier. The next step in this approach requires an algebraic closed
formula condition for the intersection of two space curves, one given
implicitly and the other one given by a parameterization.

In this paper we build resultant techniques to get original and quite
complete answers to that question and to similar ones. We think that
some of our results can lead to efficient and precise procedures of
practical interest, other ones have only a theoretical interest.

We tried to keep the paper understandable by a large audience of
readers. However the general idea developed here is more easily
expressed using homological algebra and algebraic geometry. We have
explicited the computations in important special cases, and we
provided bounds on size and degrees of the expressions.

The paper is divided in four sections organized as follows. 
Section 2 describes  the application  from  
which emerged our question. In section 3 one recalls some basic definitions 
and  the homological technique  we will use for our developments. 
We present our main new tool the determinantal resultant in section
4. In section 5, we prove several propositions and theorems which 
  provide an answer, for each case,  to our computational question;
  this last 
section ends with a simple  example on which we illustrate our different
 algorithms.

\section{Motivations : semi-implicit representation of surfaces}

An important question in Computer Geometric  Aided Design is to compute
precisely and efficiently the intersection of two patches of algebraic
surfaces given by rational parametrizations. In \cite{Gal03}, we noted
that it is much easier to represent one of the two patches by a family
of implicitly defined curves than to compute an implicitization of 
the whole surface. This technique has for instance the advantage,
among others, to 
avoid the difficulties coming from the base points when
implicitizing. By generic flatness, except for few curves, 
 this implicitization step can be done in an  uniform way for all curves.
We called that process semi-implicit representation of a parameterized surface.
Moreover in many cases of interest for the applications, these curves
have extra algebraic properties such as being determinantal.
So we are led to an apparently simpler question but in a relative 
setting (i.e. on top of a product of one dimensional spaces): detect
when two space curves intersect.

The answer to this question defines a curve in a two dimensional space, so
in theory it can be given by a generalized resultant. Providing such a
resultant is the target of this paper.

The fact that we obtain this equation via a resultant is a
 guaranty of a good 
numerical stability of the output. Moreover we can keep the result as a 
 compact straight line program made by computation of some minors.
Once we have this equation, the work for solving the CAGD problem can
be finished by (certified) numerical computations. These
numerical computations will  lift the
obtained real curve to a real curve in the four dimensional space where
live the parameters of the input parametrized surfaces.

\section{Preliminaries and previous results}

This section is devoted to some known results about multigraded
resultants. Our aim is here to recall the (co)homological techniques,
that we will heavily use in section \ref{resultant}; they were introduced in
\cite{KSZ92} to study these resultants via determinants of
complexes. As particular cases we recover the well-known Sylvester and
Dixon resultants.

\subsection{Determinants of complexes} Given a finite complex of free
$R$-modu\-les, where $R$ is a factorial regular noetherian integral
domain, which is 
generically exact (that is exact after tensorization by the field of
fractions of $R$), we can associate to it a \emph{determinant} (also
called MacRae's invariant) which measures the part of its homology
supported in codimension 1. We refer the reader to \cite{KnMu76} and
\cite{GKZ94} appendix A for the general definition and properties. In
this paper we will only need to compute determinants of complexes in
the simple case of four-term exact complexes of finite-dimensional
vector spaces. We recall a known method, going back to Cayley
\cite{Cay48}, to do this.

Suppose that we have a four-term exact complex of vector spaces
\begin{equation}\label{4term}
0\rightarrow V_3 \xrightarrow{\partial_2} V_2
\xrightarrow{\partial_1} V_1 \xrightarrow{\partial_0} V_0 \rightarrow
0.
\end{equation}
Since $\partial_0$ is surjective, $V_1$ decomposes into $V_0\oplus
V_1'$ and $\partial_1=(\begin{array}{cc} \phi_0 & m_0
  \end{array})$  with $\det(\phi_0)\neq 0$. Now, since
  $\Im(\partial_1)=\ker(\partial_0)$, $V_2$ decomposes into
  $V_1'\oplus V_3$ and $\partial_1= \left(
    \begin{array}{cc}
      m_1 & m_2 \\
      \phi_1 & m_3 \end{array}\right)$ with $\det(\phi_1)\neq
  0$. Finally since $\partial_2$ is injective and
  $\Im(\partial_2)=\ker(\partial_1)$, we have $\partial_2=\left(
    \begin{array}{c} m_4 \\ \phi_2 \end{array}\right)$ with
  $\det(\phi_2)\neq 0$. The determinant of the complex \eqref{4term}
  is then obtained as the quotient
  $\frac{\det(\phi_0)\det(\phi_2)}{\det(\phi_1)}$, and is independent
  of choices we made.
  
  Note that if $V_3=0$, we can make the same decomposition which shows
  that the determinant of \eqref{4term} is a quotient of the form
  $\frac{\det(\phi_0)}{\det(\phi_1)}$. Similarly, when moreover
  $V_2=0$, we recover the standard notion of determinant since the
  determinant of \eqref{4term} is then the determinant of the map
  $\partial_0$.

\subsection{Classical multigraded resultants}\label{multres}
We recall results about classical resultants in the
context of multigraded homogeneous polynomials, that is to say
polynomials which are homogeneous in sets of variables. Let
$l_1,\ldots,l_r$ be positive integers, and let $\xg_1,\ldots,\xg_r$ be
$r$ sets of variables consisting respectively in $l_t+1$ variables
$\xg_t=(x_{t,0},x_{t,1},\ldots,x_{t,l_t})$, for all $t=1,\ldots,r$.
We denote by $S(\dg)$, where $\dg=(d_1,\ldots,d_r)$, the vector space
of multi-homogeneous polynomials defined over the variety
$X=\PP^{l_1}\times \ldots \times \PP^{l_r}$ of multi-degree
$d_1,\ldots,d_r$, that is to say polynomials in the ring
$S=\KK[\xg_1,\ldots,\xg_r]$ which are homogeneous of degree $d_t$ in
the set of variables $\xg_t$ for all $t=1,\ldots,r$.

Let $s=\dim(X)=l_1+\ldots+l_r$, and fix $s+1$ sequences $\dg_0,
\ldots,\dg_s$, each consisting of $r$ strictly positive integers.
There exists an irreducible polynomial in
$\ZZ[S(\dg_0)\oplus\ldots\oplus S(\dg_s)]$ (unique up to sign), called
the \emph{multigraded resultant} on $X$ of type $(\dg_0,\ldots,\dg_s)$
and denoted by $\Res_{X,\dg_0,\ldots,\dg_s}$, such that, for all
$s+1-$uple $(f_0,\ldots,f_s) \in S(\dg_0)\oplus\ldots\oplus S(\dg_s)$,
$$\Res_{X,\dg_0,\ldots,\dg_s}(f_0,\ldots,f_s)=0 \Leftrightarrow
\exists x\in X \ : \ f_0(x)=\cdots=f_s(x)=0.$$
The existence of this
polynomial follows from \cite{KSZ92} section 5 (see also \cite{Stu93}
for a different presentation).  It is known to be multi-homogeneous of
degree $N_i$ in each vector space $S(\dg_i)$, $i=0,\ldots,s$. The
integers $N_i$ are determined by the intersection formula $\int_X
\prod_{j\neq i} c_1(\OO_X(\dg_j))$, where $c_1$ denotes the first
Chern class (explicit formulas are given in \cite{StZe94}).

The technique used in \cite{KSZ92} to define and study this resultant
is very powerful to derive explicit formulas, it can be summarized as
follows. The vanishing of the multigraded resultant expresses a
failure of a complex of sheaves to be exact. This allows to construct
a class of complexes of finite-dimensional vector spaces whose
determinants are the resultant. To every collection of polynomials
$f_0,f_1,\ldots,f_s$ in $S(\dg_0)\oplus\ldots\oplus S(\dg_s)$, we
associate its Koszul complex of sheaves $\Kc_\bullet$ on $X$~:
$$\Kc_{s+1}\xrightarrow{\partial_{s+1}} \Kc_{s}
\xrightarrow{\partial_s} \cdots \xrightarrow{\partial_2} \Kc_1
\xrightarrow{\partial_1} \Kc_0=\OO_X,$$
where
$\Kc_p=\wedge^p(\oplus_{i=0}^s\OO_X(-\dg_i))$, $p=0,\ldots,s+1$, and
$\partial_1=(f_0,\ldots,f_s)$.  This complex is known to be exact if
and only if $f_0,\ldots,f_s$ have no common zero in $X$; it is hence
generically exact.  To go from $\Kc_\bullet$ to complexes of
finite-dimensional vector spaces, we tensor $\Kc_\bullet$ by an
invertible sheaf $\MM=\OO_X(m_1,\ldots,m_r)$ on $X$, where
$m_1,\ldots,m_r$ are integers, then take global sections. In order to
preserve the generic exactness of $\Kc_\bullet$ when taking global
sections, the set of integers $m_1,\ldots,m_r$ have to satisfy the
following conditions.

\begin{prop}\label{det} Let $\MM=\OO_X(m_1,\ldots,m_r)$ be any
  invertible sheaf on $X$. If for all $p=0,\ldots,s+1$ and all
  integer $j > 0$ the cohomology groups
  $H^j(X,\Kc_p\otimes_{\OO_X}\MM)$ vanish, then
  $\Res_{X,\dg_0,\ldots,\dg_s}(f_0,\ldots,f_s)$ equals the determinant
  of the complex of global sections of $\Kc_\bullet\otimes_{\OO_X}
  \MM$.
\end{prop}

Let us denote by $\mg$ a sequence $(m_1,\ldots,m_r)$, and by
$K_\bullet(\mg;f_0,\ldots,f_s)$ the complex of global sections of
$\Kc_\bullet\otimes_{\OO_X} \OO_X(\mg)$. Its first differential on the
far right is naturally identified with the (generically surjective)
map $\partial_1$~:
\begin{eqnarray}\label{gen-partial}
 \oplus_{i=0}^sS(\mg-\dg_i)&  \xrightarrow{\partial_1} &
S(\mg) \\
(g_0,\ldots,g_s) & \mapsto & f_0g_0+f_1g_1+\ldots+f_sg_s. \nonumber
\end{eqnarray}
As a consequence of being the determinant of
$K_\bullet(\mg;f_0,\ldots,f_s)$, we have the corollary (see
\cite{GKZ94} theorem 34)~:

\begin{cor}\label{gcd} Under the hypothesis of proposition \ref{det}, the
  multigraded resultant $\Res_{X,\dg_0,\ldots,\dg_s}(f_0,\ldots,f_s)$
  equals the gcd of all the determinants of the square minors of size
  $\dim_\KK(S(\mg))$ of the map \eqref{gen-partial}.
\end{cor}

In fact the hypothesis needed in proposition \ref{det} (and corollary
\ref{gcd}) is purely combinatorial, for we know (see for instance
\cite{Har77} theorem III.5.1) that for all $j\geq 0$, $k\geq 1$ and
any $n\in \ZZ$, the cohomology group
\begin{equation}\label{fac}
  H^j(\PP^k,\OO_{\PP^k}(n))\neq 0 \Leftrightarrow (n\geq 0 \ 
  \mathrm{ and } \ j=0), \  \mathrm{ or } \ ( n < -k \ \mathrm{ and } \ j=k).
\end{equation}
This property gives a systematic way to deduce which integers
$m_1,\ldots,m_r$ are valid for a given multigraded resultant type. We
illustrate this technique with the two following examples of the
well-known Sylvester and Dixon resultants and refer to \cite{StZe94}
for a more systematic and detailed use of these tools to find formulas
for general multigraded resultants.

\begin{ex}\label{sylres}{\rm 
    We suppose here $r=1$ and $l_1=1$, that is $X=\PP^1$, and consider
    the resultant of two bivariate homogeneous polynomials $f_0,f_1$
    of respective degree $\dg_0:=d_1$ and $\dg_1:=d_1$ in variables
    $\xg_1=(x_0,x_1)$. Such a resultant was first studied by
    Sylvester, it corresponds to the polynomial in $\ZZ[S(d_0)\oplus
    S(d_1)]$ we have denoted $\Res_{\PP^1,d_0,d_1}(f_0,f_1)$, and is
    called the Sylvester resultant.
    
    The Koszul complex associated to $f_0$ and $f_1$ tensorized by an
    invertible sheaf $\MM=\OO_X(n)$ on $X$, with $n\in \ZZ$, is of the
    form~:
\begin{equation}\label{Sylcpx}
\OO_X(n-d_0-d_1) \xrightarrow{\binom{f_1}{-f_0}} \OO_X(n-d_0)\oplus
\OO_X(n-d_1) \xrightarrow{(f_0,f_1)} \OO_X(n). 
\end{equation}
By proposition \ref{det} we know that $\Res_{X,d_0,d_1}(f_0,f_1)$ is
the determinant of the complex of global sections of \eqref{Sylcpx},
providing $n$ is chosen such that the cohomology groups
$H^1(X,\OO_X(n-d_0))$, $H^1(X,\OO_X(n-d_1))$ and
$H^1(X,\OO_X(n-d_0-d_1))$ vanish. By \eqref{fac} we deduce easily that
$n$ has to be greater or equal to $d_0+d_1-1$, and hence that
$\Res_{X,d_0,d_1}(f_0,f_1)$ can be computed as the determinant of each
complex
\begin{equation*}
S(n-d_0-d_1) \xrightarrow{\binom{f_1}{-f_0}} S(n-d_0)\oplus S(n-d_1)
\xrightarrow{(f_0,f_1)} S(n), 
\end{equation*}
with $n\geq d_0+d_1-1$. In particular, if we choose $n=d_0+d_1-1$,
then $S(n-d_0-d_1)=\emptyset$. Thus we deduce that
$\Res_{X,d_0,d_1}(f_0,f_1)$ is the determinant of the well-known
Sylvester's matrix
$$S(d_1-1)\oplus S(d_0-1)\xrightarrow{(f_0,f_1)}S(d_0+d_1-1).$$
Notice
that $\Res_{\PP^1,d_0,d_1}$ is homogeneous in $S(d_0)$ of degree
$\int_X c_1(\OO_X(d_1))=d_1$, and homogeneous in $S(d_1)$ of degree
$d_0$.}
\end{ex}

\begin{ex}\label{dixres}{\rm 
    Here we consider the case of three bi-homo\-ge\-neous polynomials
    of same bi-degree $\dg=(d_1,d_2)$ in both sets of variables
    $\xg_1=(x_{1,0},x_{1,1})$ and $\xg_2=(x_{2,0},x_{2,1})$, defined
    on $X=\PP^1\times \PP^1$. This case was originally studied by A.L.
    Dixon in its famous paper \cite{Dix08}. The corresponding
    bi-graded resultant
    $\Res_{\PP^1\times\PP^1,\dg,\dg,\dg}(f_0,f_1,f_2)$, where
    $f_0,f_1,f_2$ are polynomials in $S(d_1,d_2)$, is called the Dixon
    resultant. To compute it we begin with the Koszul complex
    associated to these polynomials, tensorized by an invertible sheaf
    $\MM=\OO_X(m,n)$, where $m,n\in\ZZ$ ~:
    $$0 \rightarrow \OO_X(m-3d_1,n-3d_2) \rightarrow
    \OO_X(m-2d_1,n-2d_2)^3$$
    $$\rightarrow \OO_X(m-d_1,n-d_2)^3 \rightarrow \OO_X(m,n).$$
    Going
    to global sections,
    $\Res_{\PP^1\times\PP^1,\dg,\dg,\dg}(f_0,f_1,f_2)$ is then
    obtained as the determinant of this complex if the couple $(m,n)$
    is chosen such that
\begin{equation}\label{cond1}
H^j(X,\OO_X(m-id_1,n-id_2))=0 \ \  \mathrm{for \ all} \ j=1,2 \ \mathrm{and} \
i=0,1,2,3.
\end{equation}
By \emph{K\"unneth formula} (see \cite{Wei94} section 3.6) we know
that
$$H^2(X,\OO_X(m-id_1,n-id_2))= H^1(\PP^1,\OO_X(m-id_1))\otimes
H^1(\PP^1,\OO_X(n-id_2)),$$
and thus, using \eqref{fac}, we see that
condition \eqref{cond1} is satisfied when $j=2$ if either $m\geq 3d_1
-1$ or $n\geq 3d_2 -1$. Always by K\"unneth formula, the first
cohomology group $H^1(X,\OO_X(m-id_1,n-id_2))$ is the direct sum of
$$H^1(\PP^1,\OO_X(m-id_1))\otimes H^0(\PP^1,\OO_X(n-id_2))$$
and
$$H^0(\PP^1,\OO_X(m-id_1))\otimes H^1(\PP^1,\OO_X(n-id_2)).$$
Using
this equality it follows that \eqref{cond1} is satisfied for all
couples such that $(m=3d_1-1,n\geq 2d_2-1)$, $(m\geq 3d_1, n \geq 3d_2
-1)$, $(m\geq 2d_1 -1, n=3d_2 -1)$, and $(m\geq 3d_1-1, n\geq 3d_2)$.
For all these couples $\Res_{X,\dg,\dg,\dg}(f_0,f_1,f_2)$ is obtained
as the determinant of the complex
\begin{equation*}
S(m-3d_1,n-3d_2) \xrightarrow{\partial_3} S(m-2d_1,n-2d_2)^3
\xrightarrow{\partial_2} S(m-d_1,n-d_2)^3 \xrightarrow{\partial_1} S(m,n),
\end{equation*}
where the differentials  $\partial_1, \partial_2, \partial_3$ are
given by the matrices~:
$$\partial_1=\left(\begin{array}{ccc} f_0 & f_1 & f_2 \end{array}
\right), \partial_2=\left(
\begin{array}{ccc}
  0 & f_2 & f_1 \\
  f_2 & 0 & -f_0 \\
  -f_1 & -f_0 & 0
\end{array}
\right), \partial_3=\left(\begin{array}{c} f_2\\ -f_1\\ f_0
  \end{array} \right).$$
If we want to obtain $\Res_{X,\dg,\dg,\dg}(f_0,f_1,f_2)$ as the
determinant of a single matrix of type \eqref{gen-partial}, we have to
find, if it exists, a couple $(m,n)$ satisfying \eqref{cond1} and also
such that $S(m-2d_1,n-2d_2)$ is empty. It is a straightforward
computation to see that this happens for both couple $(2d_1 -1, 3d_2
-1)$ and $(3d_1 -1, 2d_2 -1)$. It follows that
$\Res_{X,\dg,\dg,\dg}(f_0,f_1,f_2)$ is the determinant of both maps
$$S(d_1-1,2d_2-1)^3 \xrightarrow{\partial_1} S(2d_1-1,3d_2-1),$$
and
$$S(2d_1-1,d_2-1)^3 \xrightarrow{\partial_1} S(3d_1-1,2d_2-1),$$
which
was given by A.L. Dixon in \cite{Dix08}, section 9.

We can also compute the multi-degree of
$\Res_{X,\dg,\dg,\dg}(f_0,f_1,f_2)$. It is homogeneous in each of the
three spaces $S(d_1,d_2)$ of degree $\int_X
c_1(\OO_X(d_1,d_2))^2=d_1d_2$, and hence of total degree $3d_1d_2$.}
\end{ex}


\section{Some determinantal resultants in low
  dimension}\label{resultant}

In this section we consider another construction of resultants which
generalizes the preceding one and was called \emph{determinantal resultants}
by the first author in \cite{Bus02} (see also \cite{BusPhD}). We first
present it, and then we further study it for the special purposes of
this article. In particular, we generalize the Sylvester
resultant to the determinantal situation that we will use in section
\ref{detect}. We also generalize the Dixon resultant to the determinantal
situation and study a determinantal resultant to detect the
intersection of two parameterized space curves. 

\subsection{Determinantal multigraded resultants}
Let $X=\PP^{l_1}\times \ldots \PP^{l_r}$, where $l_1,\ldots,l_r$ are
positive integers, and denote by $\xg_1,\ldots,\xg_r$ the 
sets of the corresponding homogeneous variables. For any sequence of
$r$ integers 
$\dg:=(d_1, \ldots,d_r)$, we denote by $S(\dg)$ the vector space of
polynomials in $S=\KK[\xg_1,\ldots,\xg_r]$ which are homogeneous of
degree $d_t$ in variables $\xg_t$ for all $t=1,\ldots,r$. We can view the
classical resultant as a tool which provides a necessary and
sufficient condition so that the matrix $(f_0,\ldots,f_s)$, defined
over $X$ of dimension $s$, is of rank lower or equal to zero in at
least one point of $X$. This interpretation leads to determinantal
resultants which give a necessary and sufficient condition so that a
given polynomial matrix, defined over $X$, is of rank lower or equal
to a given integer $p$. In what follows we focus on determinantal
resultants in the particular case where $X=\PP^{l_1}\times\ldots
\times \PP^{l_r}$, and where $p$ is the greatest possible integer.

Let $\kg_1,\ldots,\kg_n$ be $n\geq 1$ sequences of $r$ integers and
$\dg_1,\ldots,\dg_{n+s}$ be $n+s$ other sequences of $r$ integers. We
can consider the vector space $H$ of all the homogeneous map
$\oplus_{i=1}^{n+s} \OO_X(-\dg_i) \rightarrow \oplus_{i=1}^{n}
\OO_X(-\kg_i)$. It is the vector space of all the polynomial matrices
\begin{equation}\label{Hmat}
\left(\begin{array}{cccc}
h_{1,1} & h_{1,2} & \cdots & h_{1,n+s} \\
h_{2,1} & h_{2,2} & \cdots & h_{2,n+s} \\
\vdots &  \vdots &         & \vdots    \\
h_{n,1} & h_{n,2} & \cdots & h_{n,n+s} \\
\end{array}\right),
\end{equation}
where $h_{i,j} \in S(\dg_j-\kg_i)$ for all $i=1,\ldots,n$ and
$j=1,\ldots,n+s$. If all the $n(n+s)$ sequences $\dg_j-\kg_i$, with
$i=1,\ldots,n$ and $j=1,\ldots,n+s$, are sequences of strictly positive
integers, then there exists a polynomial in $\ZZ[H]$ (defined up to a
non zero constant multiple), called the determinantal resultant on
$X=\PP^{l_1}\times \ldots \times \PP^{l_r}$ of type
$(\dg_1,\ldots,\dg_{n+s};\kg_1,\ldots,\kg_n)$ and denoted by
$\Res_{X,(\dg_1,\ldots,\dg_{n+s};\kg_1,\ldots,\kg_n)}$, such that, for
all homogeneous map $\phi:\oplus_{i=1}^{n+s} \OO_X(-\dg_i) \rightarrow
\oplus_{i=1}^{n} \OO_X(-\kg_i)$,
$$\Res_{X,(\dg_1,\ldots,\dg_{n+s};\kg_1,\ldots,\kg_n)}(\phi)=0
\Leftrightarrow \exists x \in X : \rank(\phi(x)) \leq n-1.$$
This
polynomial is multi-homogeneous in the coefficients of each column of
the matrix \eqref{Hmat}, its multidegree is computed using the
Thom-Porteous intersection formula in \cite{Bus02}, proposition 2.5.

\begin{rem}\label{rem} {\rm This determinantal resultant 
  depends on the sequences $\dg_j-\kg_i$, that is to say on $H$. It is
  possible to change the $\dg_j$'s and the $\kg_i$'s without changing
  the resultant.  Note also that the case $n=1$ and
  $\kg_1=(0,\ldots,0)$ gives the classical multigraded resultant
  recalled in \ref{multres}.}
\end{rem}

A technique similar to the one exposed in \ref{multres} can be used to
compute this determinantal multigraded resultant.  It involves a
generalization of the Koszul complex, the so-called Eagon-Northcott
complex (see e.g. \cite{BrVe80}). To every homogeneous map $\phi \in H$ we
associate its Eagon-Northcott complex of sheaves $\EE_\bullet$ on $X$ 
which is of the form 
$$\EE_{s+1}\xrightarrow{\partial_{s+1}} \EE_{s}
\xrightarrow{\partial_{s}} \cdots \xrightarrow{\partial_2} \EE_1
\xrightarrow{\partial_1} \EE_0=\OO_X,$$
where
$\EE_p=\wedge^{n+p-1}E\otimes S^{p-1}F^*\otimes \wedge^{n}F^*$, for
all $p=1,\ldots,s+1$, with the notations
$E=\oplus_{i=1}^{n+s}\OO_X(-\dg_i)$ and
$F=\oplus_{i=1}^{n}\OO_X(-\kg_i)$.

\begin{prop}\label{det+} Let $\MM=\OO_X(m_1,\ldots,m_r)$ be any
  invertible sheaf on $X$. If for all $p=0,\ldots,s+1$ and all
  integer $j > 0$ the cohomology groups
  $H^j(X,\EE_p\otimes_{\OO_X}\MM)$ vanish, then
  $\Res_{X,(\dg_1,\ldots,\dg_n+s;\kg_1,\ldots\kg_n)}(\phi)$ equals the
  determinant of the complex of global sections of
  $\EE_\bullet\otimes_{\OO_X} \MM$.
\end{prop}

Let us denote by $E_\bullet(\mg;\phi)$ the complex of global sections
of the complex of sheaves $\EE_\bullet\otimes_{\OO_X} \OO_X(\mg)$,
where $\mg$ is a sequence of $r$ integers $(m_1,\ldots,m_r)$. Its
first differential on the far right is naturally identified with the
(generically surjective) map $\partial_1$~:
\begin{eqnarray}\label{gendet-partial}
 \bigoplus_{1\leq i_1<\ldots<i_n\leq
   n+s}S(\mg-\mathbf{e}_{i_1,\ldots,i_n})&  \xrightarrow{\partial_1} & 
S(\mg) \\ (\ldots,g_{i_1,\ldots,i_n},\ldots) & \mapsto & \sum_{1\leq
  i_1<\ldots<i_n\leq n+s}g_{i_1,\ldots,i_n}\Delta_{i_1,\ldots,i_n}, \nonumber
\end{eqnarray}
where $\Delta_{i_1,\ldots,i_n}$ is the determinant of the submatrix of
\eqref{Hmat} corresponding to columns $i_1,\ldots,i_n$; the sequence
$\mathbf{e}_{i_1,\ldots,i_n}$ denotes its multi-degree. We have 
the corollary ~:

\begin{cor}\label{gcd+} Under the hypothesis of proposition \ref{det+}, the
  determinantal multigraded resultant
  $\Res_{X,(\dg_1,\ldots,\dg_{n+s};\kg_1,\ldots\kg_n)}(\phi)$ equals
  the gcd of all the square minors of size $\dim_\KK(S(\mg))$ of the
  map \eqref{gendet-partial}.
\end{cor}

We now present in more detail two particular situations which
correspond to generalizations of the Sylvester resultant and of the
Dixon resultant. We will end this section with a determinantal
resultant specially designed to detect the intersection of two
parameterized space curves.

\subsection{Determinantal Sylvester resultants}\label{detsylres}
Let $n$ be a positive integer, we consider matrices with polynomial
entries in $\PP^1$ (more precisely homogeneous matrices defined in
$\PP^1$) of size $n\times (n+1)$, and give a necessary and sufficient
condition on the coefficients of these polynomial entries so that the
rank drops. The case $n=1$, that is a $1\times 2$ matrix, corresponds
to the classical Sylvester resultant. If we add one line and one
column we obtain a new Sylvester resultant corresponding to the
vanishing of all the $2\times 2$ minors of a $2\times 3$ matrix. This
is graphically illustrated by the following figure

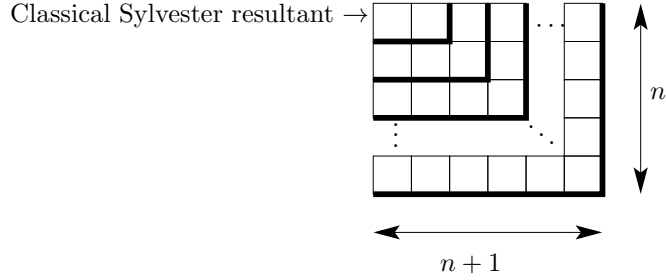
\begin{figure}[h]\label{detsylpict} 
  \vspace{.5cm} \input{sylv.pstex_t}
\end{figure}

\begin{thm} Let $n$ be a  positive integer, and 
  $(d_1,\ldots,d_{n+1})$, $(k_1,\ldots,k_n)$ be two sequences of
  integers such that $d_i-k_j>0$ for all $i,j$. For any morphism
  $\phi:\oplus_{i=1}^{n+1} \OO_{\PP^1}(-d_i) \rightarrow
  \oplus_{i=1}^{n} \OO_{\PP^1} (-k_i)$ and for any integer $$m\geq
  \sum_{j=1}^{n+1}d_j - \sum_{j=1}^{n}k_j -\min_{j=1,\ldots,n} k_j -
  1,$$
  the determinantal resultant
  $\Res_{\PP^1,(d_1,\ldots,d_{n+1};k_1,\ldots,k_n)}(\phi)$ is the
  determinant of the finite-dimensional vector spaces
  $$\bigoplus_{i=1}^n S(m+k_i-\sum_{j=1}^{n+1}d_j + \sum_{j=1}^{n}k_j)
  \xrightarrow{\phi^*} \bigoplus_{i=1}^{n+1} S(m+d_i
  -\sum_{j=1}^{n+1}d_j + \sum_{j=1}^{n}k_j) \xrightarrow{\wedge^n\phi}
  S(m).$$
  Moreover, if $k:=k_1=\ldots=k_n$ then
  $\Res_{\PP^1,(d_1,\ldots,d_{n+1};k_1,\ldots,k_n)}(\phi)$ is the
  determinant of the square matrix of size
  $d_1+\ldots+d_{n+1}-(n+1)k$,
  $$\bigoplus_{i=1}^{n+1} S(d_i-k -1 ) \xrightarrow{\wedge^n\phi}
  S\bigg(\sum_{j=1}^{n+1}(d_j - k) - 1\bigg).$$
\end{thm}
Before giving the proof of this theorem, we recall that $\wedge^n\phi$
is defined, for all integer $m$, by
$$\bigoplus_{i=1}^{n+1} \OO_{\PP^1}(d_i-\sum_{j=1}^{n+1}d_j
+\sum_{j=1}^nk_j) \xrightarrow{\wedge^n\phi} \OO_{\PP^1} \ : \ 
(g_1,\ldots,g_{n+1}) \mapsto \sum_{i=1}^{n+1}(-1)^{i-1}g_i\Delta_i,$$
where $\Delta_i$, for $i=1,\ldots,n+1$, denotes the determinant of the
matrix $\phi$ without its $\mathrm{i}^{\mathrm{th}}$ column, and that
$\phi^*$ denotes the dual of $\phi$, that is
$$\bigoplus_{i=1}^n \OO_{\PP^1}(k_i) \xrightarrow{\phi^*}
\bigoplus_{i=1}^{n+1} \OO_{\PP^1}(d_i) \ : \ (g_1,\ldots,g_n) \mapsto
(\ldots,\sum_{i=1}^ng_i\phi_{i,k},\ldots)_{k=1,\ldots,n+1}.$$

\begin{proof} Let $m$ be any integer and $\phi:\oplus_{i=1}^{n+1}
  \OO_{\PP^1}(-d_i) \rightarrow \oplus_{i=1}^{n} \OO_{\PP^1}(-k_i)$.
  The Eagon-Northcott complex of $\phi$, tensorized by the invertible
  sheaf $\OO_{\PP^1}(m)$ is of the form~:
  $$\bigoplus_{i=1}^n \OO_{\PP^1}(m+k_i-\sum_{j=1}^{n+1}d_j +
  \sum_{j=1}^{n}k_j) \rightarrow \bigoplus_{i=1}^{n+1}
  \OO_{\PP^1}(m+d_i -\sum_{j=1}^{n+1}d_j + \sum_{j=1}^{n}k_j)
  \rightarrow \OO_{\PP^1}(m).$$
  From proposition \ref{det+} and
  \eqref{fac} it comes easily that the determinant of the complex of
  global sections of this complex is the determinantal resultant of
  $\phi$ for all $m\geq \sum_{j=1}^{n+1}d_j - \sum_{j=1}^{n}k_j
  -\min_{j=1,\ldots,n} k_j - 1$, and this proves the first claim.
  Taking for instance the lowest possible value $m_0$ of $m$, we obtain
  that $\Res_{\PP^1,(d_1,\ldots,d_{n+1};k_1,\ldots,k_n)}(\phi)$ is the
  determinant of the following finite-dimensional vector spaces
  complex
  $$\bigoplus_{i=1}^n S(k_i -\min_j k_j -1) \xrightarrow{\phi^*}
  \bigoplus_{i=1}^{n+1} S(d_i-\min_j k_j -1 )
  \xrightarrow{\wedge^n\phi} S(m_0).$$
  But the vector space
  $\bigoplus_{i=1}^n S(k_i -\min_j k_j -1)$ is zero if and only if
  $k_1=\ldots=k_n$, which proves the second claim.
\end{proof}

\begin{rem} {\rm The first statement of this theorem implies that such a
  resultant is obtained either as the quotient of two determinants, or
  as the gcd of all minors of size $\dim(S(m))$ of the map $\wedge^n\phi$. \\
  Note also that, by remark \ref{rem}, we can suppose that the minimum
  of the integers $k_1,\ldots,k_n$ is zero without changing the
  resultant. For instance we can simplify the second statement of the
  theorem by taking $d_i-k$ instead of $d_i$ for all $i=1,\ldots,n+1$,
  and then $k=0$. Thus the square matrix is of size
  $d_1+\ldots+d_{n+1}$ (to be compared with the classical Sylvester
  matrix which is of size $d_1+d_2$).}
\end{rem}
 
We can also give the multidegree of the determinantal Sylvester
resultant.

\begin{prop} Let $n$ be a  positive integer, and let 
  $(d_1,\ldots,d_{n+1})$,\\ $(k_1,\ldots,k_n)$ be two sequences of
  integers such that $d_i-k_j>0$ for all $i,j$. The determinantal
  resultant $\Res_{\PP^1,(d_1,\ldots,d_{n+1};k_1,\ldots,k_n)}$ is an
  irreducible polynomial in $\ZZ[S(d_i-k_j); \forall i,j]$ which is,
  for all $i=1,\ldots,n+1$, homogeneous in $\oplus_{j=1}^nS(d_i-k_j)$
  of degree
  $$N_i:=\sum_{j=1}^{n+1}d_j-\sum_{j=1}^n k_j - d_i.$$
  It is of total
  degree
  
  $$\deg(\Res_{\PP^1,(d_1,\ldots,d_{n+1};k_1,\ldots,k_n)})=
  n\sum_{j=1}^{n+1}d_j-(n+1)\sum_{j=1}^n k_j.$$
\end{prop}

\begin{proof} We use the result of \cite{Bus02}, section 5.1, which says that
  the resultant $\Res_{\PP^1,(d_1,\ldots,d_{n+1};k_1,\ldots,k_n)}$ is
  homogeneous in the coefficients of the $\mathrm{i}^{\mathrm{th}}$
  column of the generic map $\phi:\oplus_{i=1}^{n+1} \OO_{\PP^1}(-d_i)
  \rightarrow \oplus_{i=1}^{n} \OO_{\PP^1}(-k_i)$, that is in
  $\oplus_{j=1}^nS(d_i-k_j)$, of degree $N_i$. Each integer $N_i$, for
  all $i=1,\ldots,n+1$, is obtained as the coefficient of $\alpha_i$
  of the multivariate polynomial (in variables
  $\alpha_1,\ldots,\alpha_{n+1}$) computed itself as the coefficient
  of the monomial $t^2$ in the univariate polynomial (in variable
  $t$)~:

  $$\frac{\prod_{i=1}^{n+1}(1-(d_i+\alpha_i)t)}{\prod_{i=1}^{n}(1-k_it)}.$$
  By a straightforward computation one finds $\sum_{j=1}^{n+1}d_j
  -\sum_{j=1}^n k_j - d_i$.
\end{proof}

\subsection{Determinantal Dixon resultants}

In this subsection we consider determinantal resultants in
$X=\PP^1\times \PP^1$, assuming that
$\dg_1=\dg_2=\ldots=\dg_{n+2}=(-d_1,-d_2)$, with $n$, $d_1$ and
$d_2$ positive integers, and $\kg_1=\kg_2=\ldots=\kg_n=(0,0)$. These
determinantal resultants, that we denote by
$\Res_{\PP^1\times\PP^1,(d_1,d_2)^n}$ for simplicity, can be seen as
generalizations of the well-known Dixon resultant (see example
\ref{dixres}), obtained for $n=1$~:

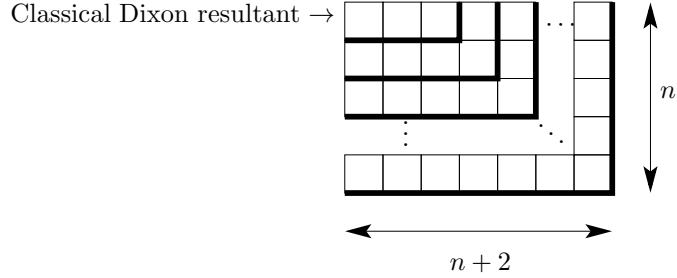
\begin{figure}[h]
  \vspace{.5cm} \input{dix.pstex_t}
\end{figure}

\begin{thm} Let $n,d_1,d_2$ be three positive integers. For any
  morphism $\phi:\oplus_{i=1}^{n+2} \OO_{\PP^1\times
    \PP^1}(-d_1,-d_2)\rightarrow \oplus_{i=1}^n\OO_{\PP^1\times
    \PP^1}$, the resultant $\Res_{\PP^1\times\PP^1,(d_1,d_2)^n}(\phi)$
  is the determinant of both square matrices
  $$\bigoplus_{i=1}^{\frac{(n+2)(n+1)}{2}}S(d_1-1,2d_2-1)
  \xrightarrow{\wedge^n\phi} S((n+1)d_1-1,(n+2)d_2-1),$$
  $$\bigoplus_{i=1}^{\frac{(n+2)(n+1)}{2}}S(2d_1-1,d_2-1)
  \xrightarrow{\wedge^n\phi} S((n+2)d_1-1,(n+1)d_2-1),$$
  of size
  $(n+2)(n+1)d_1d_2$.
  
  Moreover $\Res_{\PP^1\times\PP^1,(d_1,d_2)^n}$ is homogeneous in the
  coefficients of each column of the morphism $\phi$ of degree
  $(n+1)nd_1d_2$; its total degree is $(n+2)(n+1)nd_1d_2$.
\end{thm}

\begin{proof} Let $(p,q)$ be a couple of integers and
  $\phi:\oplus_{i=1}^{n+2} \OO_{\PP^1\times
    \PP^1}(-d_1,-d_2)\rightarrow \oplus_{i=1}^n\OO_{\PP^1\times
    \PP^1}$. The Eagon-Northcott of $\phi$, tensorized by the
  invertible sheaf $\OO_{\PP^1}(p,q)$ on $X=\PP^1\times \PP^1$, is of
  the form~:
  $$\bigoplus_{i=1}^{\frac{n(n-1)}{2}}\OO_{X}(p-(n+2)d_1,q-(n+2)d_2)
  \rightarrow \bigoplus_{i=1}^{n^2}\OO_{X}(p-(n+1)d_1,q-(n+1)d_2)$$
  $$\rightarrow
  \bigoplus_{i=1}^{\frac{(n+2)(n+1)}{2}}\OO_{X}(p-nd_1,q-nd_2)
  \xrightarrow{\wedge^n\phi} \OO_{X}(p,q).$$
  To identify the possible
  values of integers $p$ and $q$ such that the determinant of the
  complex of global sections of this complex gives the resultant
  $\Res_{\PP^1\times\PP^1,(d_1,d_2)^n}(\phi)$, we
  have to check two conditions. First $p$ and $q$ must be such that
  $$H^1(\PP^1,\OO_{\PP^1}(p-id_1))\otimes
  H^1(\PP^1,\OO_{\PP^1}(q-id_2))=0, \ \mathrm{for} \ i=n,n+1,n+2.$$
  This implies that $p\geq (n+2)d_1-1$ or $q\geq (n+2)d_2 -1$. The
  second condition is the simultaneous vanishing of
  $$H^0(\PP^1,\OO_{\PP^1}(p-id_1))\otimes
  H^1(\PP^1,\OO_{\PP^1}(q-id_2))$$
  and
  $$H^1(\PP^1,\OO_{\PP^1}(p-id_1))\otimes
  H^0(\PP^1,\OO_{\PP^1}(q-id_2))$$
  for all $i=n,n+1,n+2$. Then, the
  following sets of couples
  $$(p=(n+2)d_1-1,q\geq (n+1)d_2-1), (p\geq (n+2)d_1,q\geq
  (n+2)d_2-1),$$
  $$(p\geq (n+1)d_1-1,q= (n+2)d_2-1), \ \mathrm{and} \ (p\geq
  (n+2)d_1,q\geq (n+2)d_2),$$
  give complexes of finite-dimensional
  vector spaces whose determinant is exactly the determinantal
  resultant. In fact, among these couples of integers $(p,q)$, two of
  them reduces to a single map. They correspond to integers $p$ and
  $q$ such that the cohomology group
  $H^0(X,\OO_{X}(p-(n+1)d_1,q-(n+1)d_2))$ vanishes, that is such that
  $p\leq (n+1)d_1-1$ or $q\leq (n+1)d_2-1$. We deduce that our
  determinantal resultant is the determinant of both maps~:
  $$\bigoplus_{i=1}^{\frac{(n+2)(n+1)}{2}}S(d_1-1,2d_2-1)
  \xrightarrow{\wedge^n\phi} S((n+1)d_1-1,(n+2)d_2-1),$$
  $$\bigoplus_{i=1}^{\frac{(n+2)(n+1)}{2}}S(2d_1-1,d_2-1)
  \xrightarrow{\wedge^n\phi} S((n+2)d_1-1,(n+1)d_2-1),$$
  which give
  square matrices of size $(n+2)(n+1)d_1d_2$.
  
  The multidegree of
  $\Res_{\PP^1\times\PP^1,(d_1,d_2)^n}(\phi)$ is
  easily obtained from the determinant of both preceding maps; this
  determinantal resultant is homogeneous of degree $(n+1)nd_1d_2$ in
  the coefficients of each column of the matrix $\phi$. It is hence of
  total degree $(n+2)(n+1)nd_1d_2$.
\end{proof}

\subsection{Determinantal resultants of two parameterized space
  curves}\label{detcurveres} As we have seen in the previous
subsection, determinantal Dixon resultants deal with matrices having
bihomogeneous polynomials of bidegree $(d_1,d_2)$ entries, where
$d_1\geq 1$ and $d_2\geq 1$. The integers $d_1$ and $d_2$ are supposed
to be positive in order to fulfill a very ampleness hypothesis used to
define determinantal resultants in \cite{Bus02}, theorem 2.1. However,
this theorem states the existence of determinantal resultants in the
general setting of morphisms between two arbitrary vector bundles $E$
and $F$ on a projective and irreducible variety $X$. In some much more
simple situations, the very ampleness hypothesis (which says that the
vector bundle $\mathcal{H}om(E,F)$ is very ample on $X$) can be
weakened. We now consider such a case, which geometrically corresponds
to intersect two parameterized space curves.

Let $X=\PP^1\times \PP^1$ and consider the vector space $H$ of all the
homogeneous map $\OO_X^4 \rightarrow \OO_X(m,0)\oplus\OO_X (0,n)$,
with $m\geq 1$ and $n\geq 1$. It is the vector space of all the
polynomial matrices
\begin{equation}\label{fg}
\left(\begin{array}{cccc}
    f_0(s,t) & f_1(s,t) & f_2(s,t) &f_3(s,t) \\
    g_0(u,v) & g_1(u,v) & g_2(u,v) &g_3(u,v)
    \end{array}\right),
\end{equation}
where $f_i(s,t)$ are homogeneous polynomials of degree $m$ in $\PP^1$,
and $g_i(u,v)$ are homogeneous polynomials of degree $n$ in $\PP^1$.
We can think of both lines of these matrices as parameterized space
curves.

As we just mentioned, denoting $E=\OO_X^4$ and
$F=\OO_X(m,0)\oplus\OO_X (0,n)$, the vector bundle
$\mathcal{H}om(E,F)$ is not very ample on $X=\PP^1\times \PP^1$, but
the determinantal resultant of morphisms from $E$ to $F$ is well
defined; we will denote it by $\Res_{X,m,n}$.

\begin{thm} Let $m$ and $n$ be two positive integers. The
  determinantal resultant $\Res_{X,m,n}$ is well defined, it is
  homogeneous of degree $n$ in the coefficients of the first line, of
  degree $m$ in the coefficients of the second line, and hence of
  total degree $mn$.  For any morphism $\phi:\OO_X^4\rightarrow
  \OO_X(m,0)\oplus\OO_X (0,n)$, it is the determinant of the following
  Eagon-Northcott finite-dimensional complexes of vector spaces
  $$0\rightarrow S(p-3m;q-n)\oplus S(p-2m;q-2n)\oplus S(p-m;p-3n)
  \xrightarrow{\partial_2}$$
  $$S(p-2m;q-n)^4\oplus S(p-m;q-2n)^4 \xrightarrow{\partial_1}
  S(p-m;q-n)^6 \xrightarrow{\wedge^2\phi} S(p;q),$$
  for all $p\geq
  3m-1$ and $q\geq 3n-1$.
  
  In particular, it vanishes if and only if the rank of the $9mn\times
  24mn$ matrix
  $$S(2m-1;2n-1)^6 \xrightarrow{\wedge^2 \phi} S(3m-1,3n-1)$$
  drops.
\end{thm}

\begin{proof} First we justify the existence of $\Res_{X,m,n}$. Recall
  that we denote $E=\OO_X^4$, $F=\OO_X(m,0)\oplus \OO_X(0,n)$ and
  $H=\mathrm{Hom}(E,F)$. Although the vector bundle
  $\mathcal{H}om(E,F)$ is not very ample on $X$, the proof of
  \cite{Bus02} theorem 1 applies. In this proof the very ampleness
  hypothesis is used to show that the projection from the incidence
  variety $W$ to the projectivized parameter space $\PP(H)$
  $$W=\{(x,\phi) \in X\times \PP(H) : \rank(\phi(x))\leq 1\}
  \rightarrow \PP(H)$$
  is birational onto its image (which is called
  the resultant variety). The argument is the following~: given a
  zero-dimensional subscheme $z$ of $\PP^1\times\ \PP^1$ of degree two
  (that is two distinct points or a double point), the locus of
  matrices $\phi$ in $\PP(H)$ of rank lower or equal to 1 on $z$ is of
  codimension twice the codimension of matrices of rank lower or equal
  to 1 on only one smooth point (which is here 3). In our particular
  situation, this property remains true, even if $\mathcal{H}om(E,F)$
  is not very ample.
  
  The remaining of the proof is a standard use of techniques
  previously exposed. Choose a morphism $\phi:\OO_X^4\rightarrow
  \OO_X(m,0)\oplus\OO_X (0,n)$. Its associated Eagon-Northcott complex
  is of the form
  $$0\rightarrow \OO_X(-3m;-n)\oplus \OO_X(-2m;-2n)\oplus
  \OO_X(-m;-3n) \xrightarrow{\partial_2}$$
  $$\OO_X(-2m;-n)^4\oplus\OO_X(-m;-2n)^4 \xrightarrow{\partial_1}
  \OO_X(-m;-n)^6 \xrightarrow{\wedge^2\phi} \OO_X.$$
  After
  tensorization by an invertible sheaf $\OO_X(p;q)$, it appears, after
  a straightforward computation, that the condition of proposition
  \ref{det} to preserve generic exactness when taking global sections
  is $p\geq 3m-1$ and $q\geq 3n-1$.
  
  We can also compute the degree of this resultant using the results
  of \cite{Bus02}. Again a straightforward computation shows that
  $\Res_{X,m,n}$ is homogeneous in the coefficients of the first line
  of degree $n$, and homogeneous in the coefficients of the second
  line of degree $m$; it is hence of total degree $mn$.
\end{proof}

Note that other similar determinantal resultants exist. For instance
we can consider the determinantal resultant of $2\times 5$ matrices
with one line corresponding to a parameterized curve from $\PP^1$ to
$\PP^4$, and the second line corresponding to a parameterized surface
from $\PP^2$ (or $\PP^1\times \PP^1$) to $\PP^4$; this determinantal
resultant is thus defined over $X=\PP^1\times \PP^2$ (or, if the
second line is defined on $\PP^1\times \PP^1$, $X=\PP^1\times \PP^1
\times \PP^1$).


\section{Detecting space curves intersection}\label{detect}

Given two families of space curves depending respectively on a
parameter $\lambda$ and a parameter $\mu$, we would like to give a
necessary and sufficient condition on these parameters so that two
curves of each family intersect. Here by a parameter we mean either a
single variable or a list of independent variables. There is naturally
two ways for giving a space curve: a parametric and an implicit
formulation. In what follows we discuss the three different cases
consisting of detecting the intersection of two families of
parameterized curves, or two families of implicit curves, or finally a
parameterized curves family and an implicit curves family. We end with
an expanded example.

Note that we will always suppose hereafter that the parameterized
curves have no base points, but this is not restrictive since we can
remove base points via gcd computations.  Hereafter we denote by
$X,Y,Z,T$ the homogeneous coordinates of the projective space $\PP^3$.

\subsection{Intersection of two families of implicit curves} Let 
$\Cc_\lambda$ and $\Dc_\mu$ be two families of space curves that we
suppose here given implicitly. We would like to yield a condition on
the parameters $\lambda$ and $\mu$ so that these two families
intersect in $\PP^3$. The simplest case is to suppose that
$\Cc_\lambda$ and $\Dc_\mu$ are families of \emph{complete
  intersection} space curves. This means that $\Cc_\lambda$ (resp.
$\Dc_\mu$) is given by two polynomials, say $H_1(\lambda,X,Y,Z,T)$ and
$H_2(\lambda,X,Y,Z,T)$ (resp. $H_3(\mu,X,Y,Z,T)$ and
$H_4(\mu,X,Y,Z,T)$), homogeneous in variables $X,Y,Z,T$, which have no
common factor for all possible values of $\lambda$ (resp. $\mu$). The
condition we are looking for is then a polynomial in $\lambda$ and
$\mu$ which is nothing but the classical resultant over $\PP^3$,
$$\Res_{\PP^3}(H_1,H_2,H_3,H_4) \in \KK[\lambda,\mu],$$
which
eliminates the variables $X,Y,Z,T$. This resultant can be computed
from all the graded part $\nu$ of the first map (and the whole) Koszul
complex associated to the sequence $H_1,H_2,H_3,H_4$, if $\nu \geq
\sum_{i=1}^4\deg(H_i)-3$.
In the general situation a similar property holds, based on a known
result of elimination theory that we now recall (see \cite{Jou80} and
\cite{Laz77}). We provide a proof for the convenience of the reader,
and also to shed light on the proof of theorem \ref{bihomlazard} to
come.
\begin{thm}\label{lazard} Let $A$ be a noetherian commutative ring,
  $n\geq 1$ be a given integer, and 
  $C=A[X_1,\ldots,X_n]$, where $X_1,\ldots,X_n$ are indeterminates. We
  denote $\mm=(X_1,\ldots,X_n)$ the irrelevant ideal of $C$. Let
  $f_1,\ldots,f_r$ be $r\geq n$ homogeneous polynomials in $C$ of
  respective degree $d_1\geq d_2 \geq \ldots \geq d_r \geq 1$. Both
  following statements are equivalent~:
  \begin{itemize}
  \item[1)] $\exists \ n \in \NN$ such that $\mm^n \subset
    (f_1,\ldots,f_r)$, \\
  \item[2)] The map of free $A$-modules $\oplus_{i=1}^r C(-d_i)_\nu
    \xrightarrow{(f_1,\ldots,f_r)} C_\nu$ is surjective for all $\nu
    \geq \delta:=d_1+d_2+\ldots+d_n-n+1$.
  \end{itemize}
\end{thm}
\begin{proof} Let us denote by $K^\bullet(\fg;C)$ the Koszul complex
  associated to the sequence $(f_1,\ldots,f_r)$ in $C$. We have
  $K^q(\fg;C)=0$ for all $q>0$ and $q<-r$, and
  $K^q(\fg;C)=\wedge^{-q}(C^r)$ for $-r\leq q \leq 0$ (we here adopt
  the cohomological notation). We denote
  $H^q(\fg;C):= H^q(K^\bullet(\fg;C))$ for all $q\in \ZZ$. Note that
  since the polynomials $f_1,\ldots,f_r$ are homogeneous, all the $C$-modules
  $K^q(\fg;C)$ and $H^q(\fg;C)$ are graded. From the exact sequence of
  $A$-modules
  $$\oplus_{i=1}^r C(-d_i)_\nu \xrightarrow{(f_1,\ldots,f_r)} C_\nu
  \rightarrow H^0(\fg;C)_\nu \rightarrow 0,$$
  we deduce immediately
  that 2) implies 1), and also that we have to prove that
  $H^0(\fg;C)_\nu$ for all $\nu \geq \delta$ to show that 1) implies 2).
  
  First we denote by $C^\bullet_\mm(C)$ the standard Cech complex
  $$0\rightarrow C \rightarrow \oplus_{i=1}^nC_{X_i} \rightarrow
  \oplus_{1\leq i < j \leq n}C_{X_iX_j} \rightarrow \ldots \rightarrow
  C_{X_1\ldots X_n} \rightarrow 0,$$
  and recall that $H^i_\mm(C)\simeq
  H^i(C^\bullet_\mm(C))$. We can now construct the bicomplex
  $K^{\bullet \bullet}(\fg;C):=K^\bullet(\fg;C)\otimes
  C^\bullet_\mm(C)$ which gives two spectral sequences having the
  same limit
  $$\left\{ \begin{array}{c}
      {}^{'}E_1^{p,q}=H^q(\fg;C^p_\mm(C)) \Rightarrow E_\infty^{p,q}\\
      {}^{''}E_2^{p,q}=H^q(\fg;H^p_\mm(C)) \Rightarrow E_\infty^{p,q}
    \end{array}\right.,$$
  where, for any $C$-module $M$, $H^q(\fg;M)$ denotes
  $H^q(\fg;C)\otimes_C M$. On one hand, as we suppose 1), we deduce
  that the support of $H^\bullet(\fg;C)$ is contained in $\mm$, and
  hence that ${}^{'}E_1^{p,q}=0$ for all $p\neq 0$, and one the other
  hand ${}^{''}E_2^{p,q}=0$ for all $p\neq n$ since $H^i_\mm(C)=0$ for
  all $i\neq n$. It follows in particular that we have an isomorphism
  of graded $C$-modules $H^0(\fg;C)\simeq H^{-n}(\fg;H^n_\mm(C))$.
  Since $H^n_\mm(C)_\nu=0$ for all $\nu\geq n+1$, we deduce that
  $H^0(\fg;C)_\nu=0$ for all $\nu\geq \delta$.
\end{proof}

Let us return to our problem. We suppose now that the family $\Cc_\lambda$ is
given by $n$ polynomials $H_i(\lambda,X,Y,Z,T)_{i=1,\ldots,n}$,
homogeneous in the variables $X,Y,Z,W$ of respective degree $d_1\geq
d_2 \geq \ldots \geq d_n\geq 1$, such that for all possible values of
$\lambda$ they define a space curve. Similarly we suppose that
$\Dc_\mu$ is given by $m$ polynomials
$H_i(\mu,X,Y,Z,T)_{i=n+1,\ldots,n+m}$, homogeneous in the variables
$X,Y,Z,W$ of respective degree $d_{n+1}\geq d_{n+2} \geq \ldots \geq
d_{n+m}\geq 1$.

\begin{prop} With the above notations, let $\delta$ be the sum of the
  four greatest integers in the set $\{d_1,\ldots,d_{n+m}\}$, minus 3.
  Then for all values of $\lambda$ and $\mu$, both space curves
  $\Cc_\lambda$ and $\Dc_\mu$ intersect if and only if the map
  \begin{eqnarray*}
    \oplus_{i=1}^{n+m}\KK[X,Y,Z,T]_{\delta-d_i} & 
  \rightarrow & \KK[X,Y,Z,T]_\delta \\
   (g_1,\ldots,g_{n+m}) & \mapsto & \sum_{i=1}^{n+m}g_iH_i,
   \end{eqnarray*}
   is not surjective.
\end{prop}

The matrices involved in this proposition are in general quite big,
and almost never square. We know examine the case of two families of
parameterized curves.

\subsection{Intersection of two families of parameterized curves} Suppose given
two families $\Cc_\lambda$ and $\Dc_\mu$ of \emph{parameterized} space
curves without base points. The family $\Cc_\lambda$ corresponds to
four homogeneous polynomials of degree $m\geq 1$,
$f_{0,\lambda}(s,t),f_{1,\lambda}(s,t), f_{2,\lambda}(s,t)$,
$f_{3,\lambda}(s,t)$ without common factor for all possible value of
$\lambda$, and the family $\Dc_\mu$ to four homogeneous polynomials of
degree $n\geq 1$, $g_{0,\mu}(u,v),g_{1,\mu}(u,v), g_{2,\mu}(u,v)$,
$g_{3,\mu}(u,v)$ without common factor for all possible value of
$\mu$. We can detect their intersection with the following
proposition, using the results of subsection \ref{detcurveres}~:

\begin{prop} With the above assumptions, the determinantal resultant
  $$\Res_{\PP^1\times\PP^1,m,n}\left(\begin{array}{cccc}
      f_{0,\lambda}(s,t) & f_{1,\lambda}(s,t) & f_{2,\lambda}(s,t) &
      f_{3,\lambda}(s,t) \\
      g_{0,\mu}(u,v) & g_{1,\mu}(u,v) & g_{2,\mu}(u,v) &
      g_{3,\mu}(u,v)
    \end{array}\right)$$  
  vanishes at $\lambda_0,\mu_0 \in \KK$ if and only if
  $\Cc_{\lambda_0}$ and $\Dc_{\mu_0}$ intersect in $\PP^3$.
\end{prop}
\begin{proof} This proposition is clear since we have supposed that
  our curves have no base points. \end{proof}

However, this resultant involves big matrices.

Another way to eliminate variables $s,t,u,v$ from the six $2\times 2$
minors of the matrix
$$\left(\begin{array}{cccc} f_{0,\lambda}(s,t) & f_{1,\lambda}(s,t) &
    f_{2,\lambda}(s,t) &
    f_{3,\lambda}(s,t) \\
    g_{0,\mu}(u,v) & g_{1,\mu}(u,v) & g_{2,\mu}(u,v) & g_{3,\mu}(u,v)
    \end{array}\right)$$
  is to study, as in theorem \ref{lazard}, their associated Koszul
  complex. This is our next theorem. We  have restricted
  ourselves to the case of bi-graded polynomials of same bi-degree,
  but the technique used in the proof applies similarly to
  multi-graded polynomials of any multi-degree.
\begin{thm}\label{bihomlazard} Let $A$ be a noetherian commutative ring. Let
  $C_1=A[s,t]$ with irrelevant ideal $\mm_1=(s,t)\subset C_1$, and
  $C_2=A[u,v]$ with irrelevant ideal $\mm_2=(u,v)\subset C_2$. Define
  the bi-graded ring $C=C_1\otimes_A C_2$ with irrelevant ideal
  $\mm=\mm_1\mm_2$, and suppose given $f_1,\ldots,f_r$ 
  bi-homogeneous polynomials in $C$ of the same bi-degree $(m,n)$,
  such that $r\geq 3$, $m\geq 1$ and $n\geq 1$. 
  Both following statements are equivalent~:
  \begin{itemize}
  \item[1)] $\exists \ n \in \NN$ such that $\mm^n \subset
    (f_1,\ldots,f_r)$, \\
  \item[2)] The map of free $A$-modules $ C(-m,-n)^r_{\nu_1,\nu_2}
    \xrightarrow{(f_1,\ldots,f_r)} C_{\nu_1,\nu_2}$ is surjective for
    all $(\nu_1,\nu_2)$ such that $\nu_1\geq 3m-3$ and $\nu_2\geq 3n-3$.
  \end{itemize}
\end{thm}
\begin{proof} As in the proof of theorem
  \ref{lazard}, let us denote by $K^\bullet(\fg;C)$ the Koszul complex
  associated to the sequence $(f_1,\ldots,f_r)$ in $C$, and its
  cohomology $H^q(\fg;C):= H^q(K^\bullet(\fg;C))$ for all $q\in \ZZ$.
  Note that since the polynomials $f_1,\ldots,f_r$ are bi-homogeneous, all the
  $C$-modules $K^q(\fg;C)$ and $H^q(\fg;C)$ are naturally bi-graded.
  From the exact sequence of $A$-modules
  $$
  C(-m,-n)^r_{\nu_1,\nu_2} \xrightarrow{(f_1,\ldots,f_r)}
  C_{\nu_1,\nu_2} \rightarrow H^0(\fg;C)_{\nu_1,\nu_2} \rightarrow
  0,$$
  we have 2) implies 1), and to prove 1) implies 2) we just have
  to show that $H^0(\fg;C)_{\nu_1,\nu_2}=0$ for all $(\nu_1,\nu_2)$ 
  such that $\nu_1\geq 3m-3$ and $\nu_2\geq 3n-3$.
  
  Always as theorem \ref{lazard}, we denote by $C^\bullet_\mm(C)$ the
  standard Cech complex
  $$0\rightarrow C \rightarrow \oplus_{i=1}^4C_{m_i} \rightarrow
  \oplus_{1\leq i < j \leq 4}C_{m_im_j} \rightarrow \ldots \rightarrow
  C_{m_1m_2m_3m_4} \rightarrow 0,$$
  where $\mm=(m_1,m_2,m_3,m_4)$,
  that is $m_1=su$, $m_2=sv$, $m_3=tu$ and $m_4=tv$. Recall also that
  $H^i_\mm(C)\simeq H^i(C^\bullet_\mm(C))$.  We can now construct the
  bicomplex $K^{\bullet \bullet}(\fg;C):=K^\bullet(\fg;C)\otimes
  C^\bullet_\mm(C)$ which gives two spectral sequences having the
  same limit
  $$\left\{ \begin{array}{c}
      {}^{'}E_1^{p,q}=H^q(\fg;C^p_\mm(C)) \Rightarrow E_\infty^{p,q}\\
      {}^{''}E_2^{p,q}=H^q(\fg;H^p_\mm(C)) \Rightarrow E_\infty^{p,q}
    \end{array}\right..$$
  First, as we suppose 1), we deduce that the support of
  $H^\bullet(\fg;C)$ is contained in $\mm$, and hence that
  ${}^{'}E_1^{p,q}=0$ for all $p\neq 0$. To study ${}^{''}E_2^{p,q}$,
  we have to compute the local cohomology modules $H^i_\mm(C)$, which
  can be done (once again) with the K\"unneth formula.
  
  First, clearly $\depth_\mm(C)\geq 2$, and hence
  $H^0_\mm(C)=H^1_\mm(C)=0$. Now let $U_i=\Spec(C_i)\setminus
  V(\mm_i)$, $i=1,2$. The fibred product $U=U_1\times_{\Spec(A)} U_2$
  is nothing but $\Spec(C)\setminus V(\mm)$, and the K\"unneth formula
  gives for all $p\geq 0$~:
  $$H^p(U,\OO_U)=\bigoplus_{p_1+p_2=p} H^{p_1}(U_1,\OO_{U_1})\otimes
  H^{p_2}(U_2,\OO_{U_2}).$$
  As $H^{p+1}_\mm(C)=H^{p}(U,\OO_U)$ for all
  $p\geq 1$, and also similarly
  $H^{p+1}_{\mm_i}(C_i)=H^{p}(U_i,\OO_{U_i})$ for $i=1,2$ and $p\geq
  1$, we deduce that $H^p_\mm(C)=0$ for $p\geq 4$ and
  $$H^2_\mm(C)=C_1\otimes H^2_{\mm_2}(C_2) \oplus
  H^2_{\mm_1}(C_1)\otimes C_2, \ H^3_\mm(C)=H^2_{\mm_1}(C_1)\otimes
  H^2_{\mm_2}(C_2).$$
  Returning to our bicomplex, we deduce that
  ${}^{''}E_2^{p,q}=0$ for all $p\notin \{2,3\}$, and obtain an
  isomorphism of bi-graded $C$-modules
  $$H^0(\fg;C)\simeq H^{-2}(\fg;H^2_\mm(C))\oplus
  H^{-3}(\fg;H^3_\mm(C)).$$
  A straightforward computation using our
  knowledge of the local cohomology of the bi-graded ring $C$ shows
  that $H^0(\fg;C)_{\nu_1,\nu_2}=0$ for all $(\nu_1,\nu_2)$
  such that $\nu_1\geq 3m-3$ and $\nu_2\geq 3n-3$.
\end{proof}

\subsection{Intersection of a family of parameterized curves with a family
  of implicit curves} Let $\Cc_\lambda$ be a family of
\emph{parameterized} space curves without base points, that is four
homogeneous polynomials $f_{0,\lambda}(s,t),f_{1,\lambda}(s,t)$, 
$f_{2,\lambda}(s,t)$, $f_{3,\lambda}(s,t)$ without common factor for
all possible value of $\lambda$.

  Let $\Dc_\mu$ be another family of space curves given
  \emph{implicitly}. The simplest situation is when $\Dc_\mu$ 
  is a family of complete intersection implicit curves, that is represented by
  two homogeneous polynomials $H_{0,\mu}(X,Y,Z,T)$ and
  $H_{1,\mu}(X,Y,Z,T)$ depending on the parameter $\mu$, without
  common factor for any value of $\mu$. In this way, for any
  specialization of $\mu$ in $\KK$, the intersection of both surfaces
  defined by $H_{0,\mu}$ and $H_{1,\mu}$ is a (complete intersection)
  curve in $\PP^3$. It is then easy to see that the Sylvester
  resultant
  $$\Res_{\PP^1}\big(H_{0,\mu}(f_{0,\lambda},f_{1,\lambda},
  f_{2,\lambda},f_{3,\lambda}),H_{1,\mu}(f_{0,\lambda},f_{1,\lambda},
  f_{2,\lambda},f_{3,\lambda})\big) \in \KK[\lambda,\mu]$$
  vanishes at
  $\lambda_0,\mu_0 \in \KK$ if and only if both space curves
  $\Cc_{\lambda_0}$ and $\Dc_{\mu_0}$ intersect in $\PP^3$. This
  resultant being computed as the determinant of a square matrix (see
  \ref{sylres}), we deduce an easy algorithm to test the intersection
  of two such families of curves, but we have to suppose that the
  family $\Dc_\mu$ is a family of complete intersection curves, which
  is quite restrictive.
  
  Using the generalized Sylvester resultant we have developed in
  \ref{detsylres}, we can weaken the complete intersection hypothesis
  and obtain a similar algorithm for a wide range of families of
  implicit curves, namely the families of \emph{determinantal}
  implicit curves.

  We suppose now that the family of implicit curves $\Dc_\mu$ is a
  flat family of arithmetically Cohen-Macaulay curves in $\PP^3$. By
  the Hilbert-Burch theorem (see \cite{Eis94} theorem 20.15), such
  a family corresponds to a homogeneous map
  $$\phi_\mu:\oplus_{i=1}^{n+1}\KK[X,Y,Z,T](-d_i) \rightarrow
  \oplus_{i=1}^{n}\KK[X,Y,Z,T](-k_i),$$
  where $n$ is a positive
  integer and $d_i > k_j$ for all $i,j$, depending on the parameter
  $\mu$ such that the ideal of all its $n\times n$ minors, denoted
  $I_n(\phi_\mu)$, is of codimension 2 for all possible value of
  $\mu$. The map $\phi_\mu$ corresponds to a polynomial matrix
\begin{equation}\label{mat}
\left(\begin{array}{cccc}
\phi_{1,1,\mu} & \phi_{1,2,\mu} & \cdots &  \phi_{1,n+1,\mu} \\
\phi_{2,1,\mu} & \phi_{2,2,\mu} & \cdots &  \phi_{2,n+1,\mu} \\
\vdots & \vdots & & \vdots \\
\phi_{n,1,\mu} & \phi_{n,2,\mu} & \cdots &  \phi_{n,n+1,\mu}
\end{array}\right)
\end{equation}
where $\phi_{i,j,\mu}$ is a homogeneous polynomial in variable
$X,Y,Z,T$ of degree $d_j-k_i$, depending on the parameter $\mu$.

\begin{prop} With the above assumptions, the determinantal resultant
  
  $$\Res_{\PP^1,(d_1,\ldots,d_{n+1};k_1,\ldots,k_n)}
  \big(\phi_\mu(f_{0,\lambda},f_{1,\lambda},f_{2,\lambda},
  f_{3,\lambda})\big)$$
  vanishes at $\lambda_0,\mu_0 \in \KK$ if and
  only if $\Cc_{\lambda_0}$ and $\Dc_{\mu_0}$ intersect in $\PP^3$.
\end{prop}

Notice that, as we have seen in \ref{detsylres}, this determinantal
resultant is the quotient of two determinants. If moreover we have
$k_1=\ldots=k_n$, then it is the determinant of a single matrix.

\subsection{An example} Hereafter we apply all the
techniques developed in this section to an explicit example involving
a family of cubics and a family of conics, our aim being to compare
the different formulations we have given. We made our computations
with the software {\tt Macaulay2} \cite{M2}, using a package\footnote{available at
\url{http://math.unice.fr/~lbuse/m2package.html}} providing functions
to compute different kinds of resultants.

Our first family is a family of cubics depending on a single parameter
$\lambda$. The parametric formulation is given by
\begin{eqnarray*}
\Cc_\lambda \ : \ \PP^1 & \rightarrow & \PP^3 \\
(s:t) & \mapsto & (s^3, s^2t-t^3, \lambda s^2t+st^2, -s^3+t^3).  
  \end{eqnarray*}
  This family of cubics is in fact a determinantal family of cubics,
  their implicit equations in $\PP^3$ are obtained as the $2\times 2$
  minors of the matrix
  $$\begin{pmatrix}{X}& {X}+{Y}+{T}&
    -{X} {\lambda}-{Y} {\lambda}-{T} {\lambda}+{Z}\\
    {X}+{Y}+{T}& -{X} {\lambda}-{Y} {\lambda}-{T} {\lambda}+{Z}&
    {X}+{T}\\
      \end{pmatrix}.$$

      Our second family is a family of conics depending on a single
      parameter $\mu$. The parameterization is
\begin{eqnarray*}
\Dc_\mu \ : \ \PP^1 & \rightarrow & \PP^3 \\
(u:v) & \mapsto & (uv, \mu v^2+u^2, \mu v^2, v^2),  
  \end{eqnarray*}
  and this family is implicitly a complete intersection given by both
  equations~: 
  $$-X^2+YT-ZT, \ -\mu T+Z.$$

  Using the implicit/implicit technique to test the intersection of
  $\Cc_\lambda$ and $\Dc_\mu$ we obtain a $56\times 115$ matrix.
  Computing a $56 \times 56$ minor we get the condition we are
  looking for~:
{\small $$\begin{array}{l}  
  \lambda^{6}\mu^{3}+3 \lambda^{6}\mu^{2}-3 \lambda^{5}\mu^{3}+3
  \lambda^{6} {\mu}-6 \lambda^{5}
  \mu^{2}-\lambda^{3}\mu^{4}+\lambda^{6}-3 \lambda^{5} {\mu}-4
  \lambda^{4}\mu^{2}+4 \lambda^{3}\mu^{3}+\\
  10 \lambda^{2}\mu^{4}+ 6 {\lambda}\mu^{5}+\mu^{6}-5 \lambda^{4} {\mu}+10
  \lambda^{3}\mu^{2}+14 \lambda^{2}\mu^{3}-3 {\lambda}
  \mu^{4}-3\mu^{5}-\lambda^{4}+5 \lambda^{3} {\mu}+\\
  4 \lambda^{2}\mu^{2}-3 {\lambda}\mu^{3}+4 \mu^{4}+\lambda^{3}-2
  \lambda^{2} {\mu}+8 {\lambda}\mu^{2}-\mu^{3}+{\lambda}
  {\mu}+2\mu^{2}-\lambda+4\mu+1.\end{array}$$}

  The parametric/parametric technique is better~: we obtain a $28 \times
  48$ matrix, whose a $28\times 28$ minor gives also the desired 
  condition. Both these results are however not satisfactory because we obtain
  quite big matrices for curves of such low degrees. We now investigate the
  mixed situation parametric/implicit.
   In our case, we can consider either $\Cc_\lambda$ or $\Dc_\mu$ as the
  implicit family. We begin by considering that $\Dc_\mu$ is the
  implcit family. Our condition is then given by a classical Sylvester
  resultant~: it is obtained as the determinant of a $9 \times 9$
  matrix. Now considering that $\Cc_\lambda$ is the
  implicit family the condition we are looking for is given in a very
  compact way by the following $6\times 6$ matrix 
{\scriptsize
$$
       \left(
      \begin{array}{cccc} {-1}& 0& {\lambda}& 0\\
        -\lambda-2& {-1}& 2\lambda&
        {\lambda}\\
        -\lambda-2\mu-3& -\lambda-2& 2 {\lambda} {\mu}+3\lambda-\mu+1&
        2\lambda\\
        -{\lambda} {\mu}-\lambda-\mu-2& -\lambda-2\mu-3& 2 {\lambda}
        {\mu}+2\lambda-\mu+1&
        2 {\lambda} {\mu}+3\lambda-\mu+1\\
        -\mu^{2}-2\mu-1& -{\lambda} {\mu}-\lambda-\mu-2&
        {\lambda}\mu^{2}+2 {\lambda} {\mu}-\mu^{2}+{\lambda}-\mu&
        2 {\lambda} {\mu}+2\lambda-\mu+1\\
        0& -\mu^{2}-2\mu-1& 0& {\lambda}\mu^{2}+2 {\lambda}
        {\mu}-\mu^{2}+{\lambda}-\mu
      \end{array}
      \right.
$$
$$
\linebreak
       \left.
     \begin{array}{cc}{-\lambda^{2}}& 0\\
        -2\lambda^{2}+1&
        {-\lambda^{2}}\\
        -2\lambda^{2} {\mu}-3\lambda^{2}+2 {\lambda} {\mu}+{2}&
        -2\lambda^{2}+1\\
        -2\lambda^{2} {\mu}-2\lambda^{2}+2 {\lambda} {\mu}+{\mu}+{2}&
        -2\lambda^{2} {\mu}-3\lambda^{2}+2 {\lambda} {\mu}+{2}\\
        -\lambda^{2}\mu^{2}-2\lambda^{2} {\mu}+2
        {\lambda}\mu^{2}-\lambda^{2}+2 {\lambda}
        {\mu}-\mu^{2}+{\mu}+1&
        -2\lambda^{2} {\mu}-2\lambda^{2}+2 {\lambda} {\mu}+{\mu}+{2}\\
        0& -\lambda^{2}\mu^{2}-2\lambda^{2} {\mu}+2
        {\lambda}\mu^{2}-\lambda^{2}+2 {\lambda} {\mu}-\mu^{2}+{\mu}+1
      \end{array}
        \right)$$
}
These experimental results show, in our point of view, how relevant the
concept of semi-implicitization is in this kind of problems.
  
  \bibliographystyle{alpha}

\begin{thebibliography}{Bus01b}

\bibitem[BJ02]{BuJo02}
Laurent Bus\'e and Jean-Pierre Jouanolou.
\newblock On the closed image of a rational map and the implicitization
  problem.
\newblock {\em Preprint math.AG/0210096}, 2002.

\bibitem[Bus01a]{BusPhD}
Laurent Bus\'e.
\newblock {\em \'Etude du r\'esultant sur une vari\'et\'e alg\'ebrique}.
\newblock PhD thesis, University of Nice, 2001.

\bibitem[Bus01b]{Bus01}
Laurent Bus\'e.
\newblock Residual resultant over the projective plane and the implicitization
  problem.
\newblock {\em proceedings ISSAC2001}, pages 48--55, 2001.

\bibitem[Bus02]{Bus02}
Laurent Bus\'e.
\newblock Determinantal resultant.
\newblock {\em Preprint math.AG/0209404}, 2002.

\bibitem[BV80]{BrVe80}
Winfried Bruns and Udo Vetter.
\newblock Determinantal rings.
\newblock {\em Lecture Notes in Mathematics}, 1327, 1980.

\bibitem[Cay48]{Cay48}
A.~Cayley.
\newblock On the theory of elimination.
\newblock {\em Cambridge and Dublin Math. Journal}, 3:116--120, 1848.

\bibitem[CLO98]{CLO98}
David Cox, J.~Little, and D.~O'Shea.
\newblock {\em Using algebraic geometry}.
\newblock Graduate Texts in Mathematics. Springer, 1998.

\bibitem[Cox01]{Cox01}
David~A. Cox.
\newblock Equations of parametric curves and surfaces via syzygies.
\newblock {\em Contemporary Mathematics}, 286:1--20, 2001.

\bibitem[Dix08]{Dix08}
A.~L. Dixon.
\newblock The eliminant of three quantics in two independent variables.
\newblock {\em Proc. London Math. Soc.}, 7:49--69, 1908.

\bibitem[Eis94]{Eis94}
D.~Eisenbud.
\newblock {\em {C}ommutative {A}lgebra with a view toward {A}lgebraic
  {G}eometry}, volume 150 of {\em Graduate Texts in Math.}
\newblock Springer-Verlag, 1994.

\bibitem[Gal03]{Gal03}
Andr\'e Galligo.
\newblock Semi-implicit representation of parameterized bi-cubic surfaces and
  applications.
\newblock {\em in preparation}, 2003.

\bibitem[GKZ94]{GKZ94}
I.M. Gelfand, M.M. Kapranov, and A.V. Zelevinsky.
\newblock {\em Discriminants, {R}esultants and {M}ultidimensional
  {D}eterminants}.
\newblock Birkh{\"{a}}user, Boston-Basel-Berlin, 1994.

\bibitem[GS]{M2}
Daniel~R. Grayson and Michael~E. Stillman.
\newblock Macaulay 2, a software system for research in algebraic geometry.
\newblock Available at http://www.math.uiuc.edu/Macaulay2.

\bibitem[Har77]{Har77}
Robin Hartshorne.
\newblock {\em Algebraic Geometry}.
\newblock Springer-Verlag, 1977.

\bibitem[Jou80]{Jou80}
J.-P. Jouanolou.
\newblock Id\'eaux r\'esultants.
\newblock {\em Adv. in Math.}, 37:212--238, 1980.

\bibitem[KM76]{KnMu76}
F.~Knudsen and D.~Mumford.
\newblock The projectivity of the moduli space of stable curves. i:
  Preliminaries on {D}et and {D}iv.
\newblock {\em Math. Scand.}, 39:19--55, 1976.

\bibitem[KSZ92]{KSZ92}
M.~M. Kapranov, B.~Sturmfels, and A.V. Zelevinsky.
\newblock Chow polytopes and general resultants.
\newblock {\em Duke Math. Journal}, 67:189--218, 1992.

\bibitem[Laz77]{Laz77}
D.~Lazard.
\newblock Alg\`ebre lin\'eaire sur $k[x_1,\ldots,x_n]$ et \'elimination.
\newblock {\em Bull. Soc. math. France}, 105:165--190, 1977.

\bibitem[Stu93]{Stu93}
Bernd Sturmfels.
\newblock Sparse elimination theory.
\newblock {\em Sympos. Math., XXXIV, Cambridge Univ. Press}, pages 264--298,
  1993.

\bibitem[SZ94]{StZe94}
Bernd Sturmfels and Andrei Zelevinsky.
\newblock Multigraded resultants of {S}ylvester type.
\newblock {\em Sympos. Math., XXXIV, Cambridge Univ. Press}, 163:115--127,
  1994.

\bibitem[Wei94]{Wei94}
Charles~A. Weibel.
\newblock {\em An introduction to homological algebra}.
\newblock Cambridge studies in advanced mathematics : Cambridge University
  Press, 1994.
\end{thebibliography}

\end{document}